\documentclass[12pt,a4paper,reqno,dvipsnames]{amsart}

\usepackage{amsfonts}
\usepackage{amsmath,mathscinet}
\usepackage{amssymb,stmaryrd,amscd}
\usepackage{amsthm}
\usepackage{amsxtra}
\usepackage{bbm}
\usepackage{braket}
\usepackage{comment}
\usepackage{extarrows}
\usepackage{graphicx}
\usepackage{latexsym}
\usepackage{mathrsfs}
\usepackage{yhmath}
\usepackage{nicematrix}
\usepackage{mathtools}

\usepackage{float}
\usepackage{booktabs}
\usepackage{xcolor}
\usepackage{verbatim}
\usepackage{listings}
\usepackage[normalem]{ulem}
\usepackage{setspace}
\usepackage[pdfborder={0 0 0}]{hyperref} 
\hypersetup{colorlinks, citecolor=blue, linkcolor=blue, urlcolor=ForestGreen}
\usepackage{bm}
\usepackage{tikz}
\usepackage[initials,lite]{amsrefs} 
\AtBeginDocument{%
    \def\MR#1{}
}

\usepackage{cleveref}

\usepackage[cm]{fullpage}

\makeatother

\theoremstyle{plain}
\newtheorem{Theorem}{Theorem}[section]
\newtheorem{Lemma}[Theorem]{Lemma}
\newtheorem{Corollary}[Theorem]{Corollary}
\newtheorem{Proposition}[Theorem]{Proposition}

\theoremstyle{definition}

\newtheorem{Assumptions and Discussion}[Theorem]{Assumptions and Discussion}

\newtheorem{Example}[Theorem]{Example}
\newtheorem{Definition}[Theorem]{Definition}

\newtheorem{Remark}[Theorem]{Remark}

\theoremstyle{remark}

\newtheorem*{acknowledgment*}{Acknowledgment}

\numberwithin{equation}{section}

\usepackage[shortlabels]{enumitem}
\SetEnumerateShortLabel{a}{\textup{(\alph*)}}
\SetEnumerateShortLabel{A}{\textup{(\Alph*)}}
\SetEnumerateShortLabel{1}{\textup{(\arabic*)}}
\SetEnumerateShortLabel{i}{\textup{(\roman*)}}
\SetEnumerateShortLabel{I}{\textup{(\Roman*)}}
\def\depth{\operatorname{depth}}
\def\GrMod{\operatorname{GrMod}}
\def\gldim{\operatorname{gldim}}
\def\grade{\operatorname{grade}}

\def\Var{\operatorname{Var}}
\def\height{\operatorname{height}}
\def\gr{\operatorname{gr}}
\def\id{\operatorname{id}} 
\newcommand\Ext{\operatorname{Ext}}
\newcommand\Tor{\operatorname{Tor}}
\newcommand\Max{\operatorname{Max}}
\renewcommand\d{\partial}
\def\TT{{\mathbb T}}
\newcommand\gExt{\underline{\Ext}}
\newcommand\gHom{\underline{\Hom}}
\newcommand\oneone{(\mathbf{1},\mathbf{1})}

\def\ann{\operatorname{ann}}

\def\bar#1{\overline{#1}}

\def\CC{{\mathbb C}}

\def\conv{\operatorname{conv}}

\def\deg{\operatorname{deg}}
\def\depth{\operatorname{depth}}
\def\GrMod{\operatorname{GrMod}}

\def\dim{\operatorname{dim}}

\def\grade{\operatorname{grade}}

\def\gr{\operatorname{gr}}

\def\id{\operatorname{id}}

\def\ini{\operatorname{in}} 

\def\Mod{\operatorname{Mod}}

\def\NN{{\mathbb N}}

\def\part{\operatorname{part}}

\def\PP{{\mathbb P}}

\def\QQ{{\mathbb Q}}

\def\ZZ{{\mathbb Z}}

\newcommand\bfa{\mathbf{a}}

\newcommand\bfF{\mathbf{F}}

\newcommand\bfm{\mathbf{m}}
\newcommand\bfn{\mathbf{n}}

\newcommand\fd{\operatorname{fd}}

\newcommand\frakm{\mathfrak{m}}

\newcommand\frakp{\mathfrak{p}}

\newcommand\Hom{\operatorname{Hom}}

\newcommand\Spec{\operatorname{Spec}}

\newcommand\Supp{\operatorname{Supp}}

\newcommand{\pd}{\operatorname{pd}}

\newcommand{\h}{{(\operatorname{h})}}

\begin{document}

\title{On Cohen--Macaulay modules over the Weyl algebra}

\author[Kuei-Nuan Lin, Jen-Chieh Hsiao]{Kuei-Nuan Lin, Jen-Chieh Hsiao}




\address{Department of Mathematics, The Penn State University, McKeesport, PA,
15132, USA}
\email{linkn@psu.edu}

\address{Department of Mathematics, National Cheng Kung University, Tainan City, Taiwan}
\email{jhsiao@mail.ncku.edu.tw}

\begin{abstract}
We propose a definition of Cohen--Macaulay modules over the Weyl algebra $D$ and give a sufficient condition for a GKZ $A$-hypergeometric $D$-module to be Cohen--Macaulay. 
\end{abstract}

\maketitle

\section{Introduction}
The main purpose of this article is to provide a well-defined and reasonable definition of Cohen--Macaulay $D$-modules, where $D=k[x]\langle \d \rangle = k[x_1, \dots, x_n]\langle \d_1, \dots, \d_n \rangle $  is the $n^{\mathrm{th}}$ Weyl algebra over a field $k$ of characteristic~$0$. 
The concept of Cohen--Macaulayness is essential in commutative algebra \cites{Eisenbud_book,BH93}, combinatorics \cites{Stanley_book,MS05_book}, and algebraic geometry \cites{GW20_book,KM98_book}. For noncommutative rings, there are different notions of Cohen--Macaulayness in the literature, see for instance \cite{Rump} and the references therein. Our definition of Cohen--Macaulay $D$-modules is closer to the properties studied in \cites{SZ94,Lev92} and is compatible with the commutative algebraic definition under Gr\"obner deformation \cite{SST_book}. Notably, it allows for computer-based verification of the Cohen--Macaulay condition.

Readers unfamiliar with $D$-module theory are referred to \cites{Bjork,SST_book} for an introduction and its significance.
One of the most influential topics in this theory is the GKZ $A$-hypergeometric systems, in which the Cohen--Macaulay condition also plays a crucial role \cites{GKZ89,MMW05}. As substantial examples, we give a sufficient condition for a GKZ system to be Cohen--Macaulay as a $D$-module.

Let us outline the main ideas as follows.
\subsection{} We first recall  the Cohen--Macaulay condition in the commutative setting;
see Sections~\ref{sec:homalg},\ref{sec:commCM} for more details.

For a commutative Noetherian local ring $R$, a nonzero finite $R$-module $M$ is {\it Cohen--Macaulay} if  $\depth(M) =\dim(M)$. We call $R$ a  Cohen--Macaulay ring if it is Cohen--Macaulay $R$-module over itself. The depth of $M$ is the maximal length of any regular sequence on $M$, which is related to the projective dimension of $M$ by the Auslander--Buchsbaum formula:
\begin{equation}\label{eq:formula1}
    \pd(M)+\depth(M) = \depth(R).
\end{equation}
On the other hand, the dimension of $M$
is related to the grade $j(M)$ of $M$ by
\begin{equation}\label{eq:formula2}
    j(M)+\dim(M) = \dim(R).
\end{equation}
 Therefore, when $R$ is Cohen--Macaulay, a nonzero finite $R$-module $M$ is Cohen--Macaulay if and only if $j(M)=\pd(M)$. 

In general, we say that a nonzero finite module $M$ over a commutative Noetherian ring $R$ is {\it Cohen--Macaulay} if the localization $M_\frakp$ is a Cohen--Macaulay $R_\frakp$-module for every prime ideal $\frakp$ in the support $\Supp(M)$ of $M$. 

In the case where $R=\bigoplus_{i\ge 0} R_i$ is an $\NN$-graded algebra over a field $k=R_0$ with maximal graded ideal $\frakm =\bigoplus_{i\ge 1} R_i $, there is a parallel theory of graded Cohen--Macaulay modules as in the local case. Formulas analogous to \eqref{eq:formula1} and \eqref{eq:formula2} hold for nonzero finite graded $R$-modules. A nonzero finite graded $R$-module $M$ is Cohen--Macaulay if and only if $M_\frakm$ is Cohen--Macaulay over the local ring $R_\frakm$.

\subsection{}
The reasons behind the definition of Cohen--Macaulay $D$-modules are discussed in Sections \ref{sec:resD},\ref{sec:CMDmod}. We summarize them here.

Consider the homogenized Weyl algebra $D^\h$. Filter $D$ and $D^\h$ by Bernstein filtration, namely the $\oneone$-filtration, and denote the associated graded algebras by $S=\gr D$ and $S[h]= \gr D^\h$. Both $S$ and $S[h]$ are polynomial algebras over $k$. We consider the three $k$-algebra $S$, $S[h]$, and $D^\h$ as $\NN$-graded algebras defined by the total degree of their elements. For a $D$-module $M$ or a graded $D^\h$-module $M^\h$, the associated graded module obtained from the $\oneone$-filtration on $D$ or $D^\h$ is denoted by $\gr M$ or $\gr M^\h$, respectively.

Let $M$ be a nonzero finite $D$-module and $\bfF$ a free resolution of $M$ over $D$ that induces a graded minimal free resolution $\gr \bfF$ of $\gr M$ over $S$. By results in \cite{OT01}, there exists a graded $D^\h$-module $M^\h$ that admits a free resolution $\bfF^\h$ of $M^\h$ satisfying the following properties.
\begin{itemize}
    \item The dehomogenization of $\bfF^\h$ coincides with $\bfF$.
    \item The resolution $\bfF^\h$ induces a graded minimal free resolution $\gr \bfF^\h$ of $\gr M^\h$ over $S[h]$.
\end{itemize}
It is observed in Proposition~\ref{prop:GB_Betti} that $\gr \bfF^\h = k[h] \otimes \gr \bfF$ and that the graded Betti numbers of the three modules $M^\h$, $\gr M$, and $\gr M^\h$ are identical. In particular, these modules have the same projective dimension. Since $\gr M$ and $\gr M^\h$ are graded modules over polynomial rings, the formulas analogous to \eqref{eq:formula1} and \eqref{eq:formula2} hold for them. In particular, the $S$-module $\gr M$ is Cohen--Macaulay if and only if  $\gr M^\h$ is Cohen--Macaulay over $S[h]$.  

On the other hand, for the graded $D^\h$-module $M^\h$, a formula analogous to \eqref{eq:formula2} can be obtained by analyzing the structure of $M^\h$ as a module over the filtered ring $D^\h$ (See Theorem~\ref{thm:d+j=w} or \cite{Bjork}). Moreover, the $\NN$-graded structure on $D^\h$ guarantees an  Auslander–Buchsbaum formula in this noncommutative setting (See Theorem~\ref{thm:noncomAB} or \cite{Jor98}). Consequently, the condition $j(M^\h)= \pd(M^\h)$ is equivalent to the condition of $\gr(M^\h)$ being Cohen--Macaulay. Based on these observations, we make the following definition (Definition~\ref{def:CMD}):
\begin{center}
    A finite left $D$-module $M$ is {\it Cohen--Macaulay} if  $\gr M$  is a Cohen--Macaulay $\gr D$-module.
\end{center}

\subsection{}
To obtain examples, we are interested in the characterization of a GKZ $A$-hypergeometric $D$-modules $M_A(\beta) = D/H_A(\beta)$ to be Cohen--Macaulay, where $H_A(\beta)$ is the left $D$-ideal generated by the toric ideal $I_A \subset \CC[\d]=\CC[\d_1, \dots,\d_n]$ and the Euler operators $E-\beta$. We show in Theorem~\ref{thm:CM_GKZ} that:
\begin{center}
    if the $\CC[\d]$-module $\CC[\d]/\ini I_A$ is Cohen--Macaulay, then $M_A(\beta)$ is a Cohen--Macaulay $D$-module.
\end{center}
Here, the initial ideal $\ini I_A$ of $I_A$ is obtained by the standard grading on the polynomial ring $\CC[\d]$. The following two steps achieve this theorem.
\begin{itemize}
    \item Apply the results in \cite{SW08} concerning the associated prime of $\ini I_A$ to show that $\ini_{\oneone}(E-\beta)$ form a linear system of parameters on $\CC(x)[\d]/I_A$.
    \item Extend a argument in \cite{SST_book} to show that $\ini_{\oneone} (E-\beta)$ form a regular sequence on $\gr(D/DI_A)$. This implies that $\gr M_A(\beta) = S/\ini_{\oneone}H_A(\beta)$ is Cohen--Macaulay.
\end{itemize} 

We note that the Cohen--Macaulay condition on the toric ring $\CC[\d]/I_A$ is equivalent to the absence of rank-jumps in the GKZ hypergeometric system $H_A(\beta)$ \cite[Corollary~9.2]{MMW05}. It is known that if the Gr\"obner deformation $\CC[\d]/\ini I_A$ is Cohen--Macaulay, then so is $\CC[\d]/I_A$ \cite[Proposition~1.6.2]{BCRV}. However, the converse of this statement is not true in general (\Cref{CMExample}). We also note that there exists a toric ideal $I_A$ with $\CC[\d]/ \ini I_A$ not Cohen--Macaulay, but the corresponding GKZ $D$-module $M_A(\mathbf{0})$ is Cohen--Macaulay (\Cref{CMExample}).

\section{Preliminaries on homological algebra} \label{sec:homalg}
In this section, rings are always associative with unit elements and modules (left or right) are unitary. 
\subsection{} We recall the definitions and some standard facts about the projective and global dimensions of modules. See \cite[Chapter 4]{Wei94} for more details. 

Let $A$ be a ring. We work in the category $\Mod(A)$ of left $A$-modules. The {\it projective dimension} $\pd(M)$ of an $A$-module $M$ is the minimum integer $n$ (if it exists) such that there is a resolution of $M$ by projective modules
$ 0 \rightarrow P_n \rightarrow \cdots \rightarrow P_1 \rightarrow P_0\rightarrow M \rightarrow 0.$ Equivalently, we have \[\pd(M)= \sup\{ i : \Ext^i_A(M,N) \ne 0 \text{ for some  }N \in \Mod(A)\}.\]
Similarly, one defines the injective dimension $\id(M)$ and the flat dimension $\fd(M)$ of $M$.

The {\it global dimension} $\gldim(A)$ is defined as the supremum of $\pd(M)$ for any $A$-module $M$. It is equal to any of the following numbers.
\begin{align*}
    \gldim(A):&= \sup\{ \pd(M) : M \in \Mod(A)\}\\
    & =\sup\{ \id(M) : M \in \Mod(A)\} \\
    &=\sup\{ \pd(A/I) : I \text{ is a left ideal of }A\} \\
    &=\sup\{ i : \Ext^i_A(M,N) \ne 0 \text{ for some }M,N \in \Mod(A)\}.
\end{align*} 
In particular, the global dimension $\gldim(A)$ can be computed by finitely generated $A$-modules. As a fundamental example, Hilbert syzygy theorem states that the polynomial ring $k[x_1, \dots, x_n]$ over a field $k$ has global dimension $n$. 

We will only consider rings that are left and right Noetherian. In this case, the flat dimension $\fd(M)$ is equal to the projective dimension $\pd(M)$ for any finitely generated (left or right) $A$-module $M$. By the Tor dimension theorem, the global dimension $$\gldim(A) = \sup\{ i : \Tor_i^A(M,N) \ne 0 \text{ for some }M,N \in \Mod(A)\}$$ is equal to the global dimension of $A$
 in the category of right $A$-modules. Another important example relevant to us is a theorem of Roos, which states that the $n^\mathrm{th}$ Weyl algebra $k\langle x_1, \dots, x_n, \d_1, \dots, \d_n \rangle$ over a field $k$ of characteristic $0$ has global dimension $n$ \cite[Chapter 2, Theorem~3.15]{Bjork}.

\subsection{} We collect some facts for modules over filtered rings. The main reference of this section is \cite[Chapter 2]{Bjork}.

Consider a ring $A$ equipped with a filtration $\Sigma_0 \subset \Sigma_1 \subset \cdots$ such that $\Sigma_0$ contains the unit elements of $A$, the ring $A = \bigcup_i \Sigma_i$, and $\Sigma_i \Sigma_j \subset \Sigma_{i+j}$ for all $i,j$. We assume that the associate graded ring $\gr(A) = \Sigma_0 \oplus \Sigma_1/\Sigma_0 \oplus \Sigma_2/\Sigma_1 \oplus \cdots$ is commutative, Noetherian, and regular of pure dimension~$\omega$.  Let $\mu = \gldim(A)$ be the global dimension of $A$. Since $\gr(A)$ is regular, the global dimension $\gldim(\gr(A)) = \omega$. By \cite[Chapter 2, Theorem~3.7]{Bjork}, we have 
\begin{equation*}
    \mu =\gldim(A) \le  \gldim(\gr(A)) =\omega.
\end{equation*}

For a finite left $A$-module $M$, its dimension is defined as $$d(M):= \sup \{ d_\frakm(M) : \frakm \in \Max(A) \}$$ where $\Max(A)$ is the set of maximal ideals of $\gr(A)$ and $d_\frakm(M) := d (\gr_\Gamma(M)_\frakm)$ is the local dimension of the associated graded module $\gr_\Gamma(M)$ of $M$ with respect to any good filtration $\Gamma$ of $M$. 
We also consider the {\it grade}  of $M$: \[\grade(M)=j(M):= \inf\{ i : \Ext^i_A(M,A) \ne 0 \}.\]
These numbers are related by 
\begin{equation*} 
0 \le j(M) \le \pd(M) \le \mu \le \omega.
\end{equation*}
Moreover, we have the following theorem \cite[Chapter 2, Theorem~7.1]{Bjork}.
\begin{Theorem}\label{thm:d+j=w}  For all nonzero finite left $A$-module $M$, we have
$d(M)+j(M) = \omega$.
\end{Theorem}
\begin{Remark}\label{rem:1}
    \begin{enumerate}
        
        \item In the special case where $A = \Gamma_0 = \gr(A)$, we have $\omega = \mu$ and the equality $d(M)+j(M) = \omega$ still holds true, where $d(M) = \dim(M) = \dim (A/\ann(M))$ is the dimension of the support of $M$.

        \item When the filtered ring $A$ is also a $\NN$-graded algebra over a field $k$, the dimension $d(M)$ is equal to the Gelfand-Kirillov dimension of $M$. An $\NN$-graded $k$-algebra satisfying Theorem~\ref{thm:d+j=w} is called {\it Cohen--Macaulay} in \cites{Lev92,SZ94}.
        \item As a consequence of Theorem \ref{thm:d+j=w}, we have the Bernstein inequality
        \[\omega- \mu \le d(M) \le \omega.\] We say that a finitely generated left $A$-module $M$ is  {\it holonomic} if $d(M) = \omega-\mu$. Equivalently, a finitely generated left $A$-module $M$ is holonomic if $i=\mu=j(M)$ is the only index such that $\Ext^i_A(M,A) \ne 0$.
    \end{enumerate}
\end{Remark}
We need the following result later, see \cite[Chapter 2, Corollary~7.5]{Bjork}.
\begin{Corollary}\label{cor:1} For a nonzero finite left $A$-module $M$, we have
    $\Ext_A^v(\Ext_A^i(M,A),A)= 0 $ if $v <i$. In particular, the grade $j (\Ext^i(M,A)) \ge i$ for any $i \ge 0$.
\end{Corollary}
\begin{Remark}\label{rem:A-regular}
    For any nonzero $A$-submodule $N$ of $\Ext^i(M,A)$ we have $d(N) \le d(\Ext^i(M,A))$, so
    it follows from Theorem~\ref{thm:d+j=w} and Corollary~\ref{cor:1} that  $j(N) \ge i$. Therefore, the algebra $A$ is Auslander-regular in the sense of \cite[Definition~2.1]{Lev92}.
\end{Remark}

\section{Cohen--Macaulay properties in commutative algebra}\label{sec:commCM}
In this section, we recall the definitions and some basic facts about Cohen--Macaulay modules over commutative Noetherian rings. See \cite{BH93} for more details.

\subsection{} Let $(R,\frakm,k)$ be a Noetherian local ring and $M$ a nonzero finite $R$-module. The {\it depth} of $M$ is defined as 
\[\depth(M):= \inf\{i : \Ext^i_R(k, M) \ne 0\},\] which is the length of any maximal $M$-sequence in $\frakm$. In particular, we have $\depth(M) \le \dim(M)$. We call $M$ a {\it Cohen--Macaulay} $R$-module if $\depth(M) = \dim(M)$. A theorem of Auslander--Buchsbaum \cite[Theorem~1.3.3]{BH93} states that if $\pd(M)<\infty$, then
\[\pd(M)+\depth(M) = \depth(R).\]
If $R$ is regular or if, more generally, $R$ is Cohen--Macaulay as a module over itself, we have $\depth(R)= \dim(R)=\gldim(R)$. In this case, it follows from Remark~\ref{rem:1}(1) that $M$ is Cohen--Macaulay if and only if $j(M)= \pd(M)$. 

In general, we say that a nonzero finite module $M$ over a Noetherian ring $S$ is {\it Cohen--Macaulay}, if $M_\frakp$ is a Cohen--Macaulay $S_\frakp$ for every prime ideal $\frakp$ in the support $\Supp(M)$ of $M$.

\subsection{} \label{sec:comm_graded}
Consider now a Noetherian $\NN$-graded ring $(R=\bigoplus_{i\ge 0} R_i, \frakm, k)$, where $k=R_0$ is a field and  $\frakm= \bigoplus_{i\ge 1} R_i$ is the maximal graded ideal of $R$.
We work in the category $\GrMod(R)$ of graded $R$-modules. For $M,N \in \GrMod(R)$, denote by  $\Hom(M,N)$ the group of all graded homomorphisms of degree $0$ between $M$ and $N$. Consider also the graded homomorphism group 
\[\gHom(M,N) = \bigoplus_{i \in \ZZ} \Hom(M,N(i)).\] When $M$ is finite, the group $\gHom(M,N)$ is equal to the ungraded homomorphism group $\Hom_R(M,N)$. Similarly, we consider the graded Ext groups $\gExt^i(M,N)$ and use them to make the following definitions.
\begin{Definition}\label{def:graded pd,j,depth} Let $M$ be a nonzero finite graded $R$-module, define    \begin{align*}
    \pd(M) &:=   \sup\{ i : \gExt^i(M,N) \ne 0 \text{ for some  }N \in \GrMod(R)\},  \\
    j(M) &:=  \inf\{ i : \gExt^i(M,R) \ne 0 \},   \\
    \depth(M) &:= \inf\{i : \gExt^i(k, M) \ne 0\}.
\end{align*}
\end{Definition}
\noindent We have $\dim(M)= \dim(M_\frakm)$, $\pd(M)= \pd(M_\frakm)$, $j(M)= j(M_\frakm)$, and $\depth(M)= \depth(M_\frakm)$.

\begin{Remark}
    When $R=k[x_1, \dots, x_n]$ is the polynomial algebra over $k$ graded by the total degree of its elements, we still have $\depth(R)= \dim(R)=\gldim(R)=n$. In this case, the Auslander--Buchsbaum formula $\pd(M)+\depth(M) = \depth(R)$ also holds for any finite graded $R$-module $M$. Moreover, by \cite[Exercise~2.1.27]{BH93} a finite graded $R$-module $M$ is Cohen--Macaulay if and only if $\depth(M) = \dim(M)$ or equivalently $j(M)=\pd(M)$. We also remark that an analogous formula of Remark~\ref{rem:1}(1) \[\dim(M) +j(M) = \gldim(R)\] also holds in the graded category $\GrMod(R)$. Therefore, for a finite graded $R$-module $M$, the notions of grade $j(M)$ defined in the two categories $\GrMod(R)$ and $\Mod(R)$ coincide.
\end{Remark}

\subsection{}
Let $(R,\frakm, k) $ be either the local ring or the $\NN$-graded ring discussed above.
Let $M$ be a finite (graded) $R$-module whose dimension $\dim M = d$. A theorem of Grothendieck states that if the $i^\mathrm{th}$ local cohomology $H^i_\frakm(M)\ne 0$ then $\depth M \le i \le d$. As a consequence, the $R$-module $M$ is Cohen--Macaulay if and only if the only nonvanishing local cohomology of $M$ supported at $\frakm$ is $H^d_\frakm(M)$ \cite[Theorem~3.5.7, Remark~3.6.18]{BH93}.

\section{Auslander--Buchsbaum formula in noncommutative graded case}\label{sec:noncom AB}
Let us state the main result in \cite{Jor98} for the convenience of the readers.

Let $A = \bigoplus_{i \in \ZZ} A_i$ be a non-commutative $\NN$-graded left Noetherian $k$-algebra, where $k=A_0$ is a field. We work in a parallel setting as in the commutative graded case (see section~\ref{sec:comm_graded}). Let $\GrMod(A)$ be the category of left graded $A$-module. For $M,N \in \GrMod(A)$, consider the graded $\gHom(M,N)$ and $\gExt^i(M,N)$. For a nonzero finite graded left $A$-module $M$, we define as in Definition~\ref{def:graded pd,j,depth} the projective dimension, $\pd(M)$, the grade, $j(M)$, and the depth, $\depth(M)$, of $M$.

Consider $A$ as a left module over itself. We assume that $\gExt^i(k,A)$ is finite dimensional over $k$ for each $i\le \depth(A)$. This condition is called the $\chi^\circ$-condition in \cite{Jor98}, which is vacuous if $\depth(A)=\infty$. Then we have the following theorem \cite[Theorem~3.2]{Jor98}.

\begin{Theorem}\label{thm:noncomAB}
    For any nonzero finite graded left $A$-module $M$ with $\pd(M) < \infty$, we have
    \[\pd(M) +\depth(M) = \depth(A).\]
    In particular, if $\gldim(A) < \infty$, by applying the above theorem to $M=k$ we get $\gldim(A) = \depth(A)$.
\end{Theorem}

We mention that there is also a local cohomology theory in this noncommutative graded setting \cite{Jor97}.

\section{Minimal free resolutions of $D$-modules}\label{sec:resD}
We summarize some results in \cite{OT01} and deduce Proposition~\ref{prop:GB_Betti} for later use.

Let $D=k[x_1, \dots,x_n]\langle \partial_1, \dots, \partial_n \rangle$ be the $n^\mathrm{th}$ Weyl algebra over a field $k$ of characteristic $0$, which is an associative  $k$-algebra generated by $x_1, \dots, x_n$ and $\partial_1, \dots, \partial_n$ with relations \[x_i x_j = x_j x_i,\quad \partial_i \partial_j = \partial_j \partial_i,\quad \partial_i x_j- x_j\partial_i= \delta_{ij}\] for $i,j= 1, \dots, n$. We will sometimes abbreviate $D=k[x]\langle \partial \rangle$ when convenient.
Consider also the homogenized Weyl algebra $D^\h = k[h,x]\langle \partial \rangle$ with relations, for $1\leq i,j\leq n$,
\[x_i x_j = x_j x_i,\quad \partial_i \partial_j = \partial_j \partial_i,\quad \partial_i x_j- x_j\partial_i= \delta_{ij}h^2,\quad hx_i=x_ih,\quad h\d_i=\d_ih.\]

\subsection{}
For $u, v, \alpha, \beta \in \ZZ^n$, a weight vector $L=(u,v)$ is called admissible for $D$ and $D^\h$ if $u_i+v_i \ge 0$. Fix any admissible weight vector $L$. The $L$-degree of a monomial $x^\alpha \d^\beta=x_1^{\alpha_1} \cdots x_n^{\alpha_n}\d_1^{\beta_1}\cdots \d_n^{\beta_n} \in D$ as well as a monomial 
$h^l x^\alpha \d^\beta \in D^\h$, $l \in \NN$, is defined by $L\cdot (\alpha, \beta)=\sum_{i} u_i\alpha_i+\sum_i v_i\beta_i $. For an element $P \in D$ or $D^\h$, its $L$-degree $\deg^L(P)$ is defined as the maximum degree of its monomials. This 
induces an increasing filtration $\{L_i\}_{i \in \ZZ}$ on $D$ and $D^\h$ by $L_iD = \{P \in D \mid \deg^L(P) \le i\}$ and $L_iD^\h = \{P \in D^\h \mid \deg^L(P) \le i\}$, respectively. For wight vectors $L$ satisfying $u_i +v_i >0$, $i= 1, \dots, n$, the associated graded rings 
\[ S:= \gr^L D = \bigoplus_{i \in \ZZ} L_i D /L_{i-1}D \cong k[x,\d], \quad S[h]:= \gr^L D^\h = \bigoplus_{k \in \ZZ} L_i D^\h /L_{i-1}D^\h \cong k[h,x,\d]\] are polynomial rings, whereby abuse of notation we still use $\d_i$ as the image of $\d_i$ in $\gr^L D$ or $\gr^LD^\h$.

For a left ideal $I$ of $D$, the weight vector $L$ induces a filtration $\{I \cap L_iD \}_{k \in \ZZ}$ on $I$. The associated graded module $\gr^L I$ can be identified as the ideal $\ini_{(u,v)} I:=\langle \ini_{(u,v)} P \mid P \in I \rangle$ in $k[x,\d]$ where $\ini_{(u,v)}P$ is the $L$-initial form of $P$. On the other hand, the $L$-filtration $L_i(D/I) = \{P +I \mid P \in L_iD\}$, ${i \in \ZZ}$, on the $D$-module $D/I$ gives the associated graded module 
\[ \gr^L(D/I) = \gr^LD /\gr^L I = S/ \ini_{(u,v)} I.\] Similarly, for a left ideal $I$ of $D^\h$ we have \[ \gr^L(D^\h/I)  = S[h]/ \ini_{(u,v)} I.\] 
\subsection{}
More generally, the $L$-filtration induces many good filtrations on modules of finite type over $D$ or $D^\h$ as follows. We only present the case of $D$-modules. The case of $D^\h$-modules can be treated similarly.

For a free $D$-module $D^r$ and for $\bfm=(m_1, \dots,m_r) \in \ZZ^r$, any admissible weight vector $L$ induces a filtration on $D^r$ by \[L_i[m]D^r = \{P \in D^r \mid \deg^L(P_j)+m_j \le i, \forall j \},\, i\in \ZZ.\] The associated graded module, denoted by $\gr^L[\bfm]D^r$,  with respect to this filtration, is a free module over $S= \gr^L D$. For a $D$-module $M$ admitting a surjective $D$-module homomorphism $D^r \xrightarrow{\varphi} M$, this induces a filtration on $M$ given by $\varphi ( L_i[m]D^r)$, $i\in \ZZ$, whose associated graded module is denoted by $\gr^L[m](M)$.
\vspace*{-0.7cm}

\subsection{}
Let us consider only the $L=\oneone$ filtration on $D$ and $D^\h$. We consider the standard $\NN$-grading on the polynomial rings $S$ and $S[h]$. The homogenized $D^\h$ is also an $\NN$-graded $k$-algebra with grading defined by the total degree of its elements.
Let $M$ be a finite $D$-module. By \cite[Proposition~4.1]{OT01}, there exists a free resolution of $D$-modules
    \[\mathbf{F}:\quad \cdots \xrightarrow[]{\varphi_2} D^{r_1} \xrightarrow[]{\varphi_1} D^{r_0} \xrightarrow{\varphi_0} M \rightarrow 0\]
    and $\bfn_i \in \ZZ^{r_i}$, such that, for $i >0$ and $j\in \ZZ$, we have
    \[\varphi_i \left( L_j [\bfn_i]D^{r_i} \right) \subseteq  L_j [\bfn_{i-1}]D^{r_{i-1}},\]
and   the induced   free resolution
\begin{align*}
\gr^L\mathbf{F}:\quad \cdots \xrightarrow[]{\bar{\varphi_2}} \gr^L[\bfn_1]D^{r_1}\xrightarrow[]{\bar{\varphi_1}} \gr^L[\bfn_0]D^{r_0} \xrightarrow[]{\bar{\varphi_0}} \gr^L[\bfn_0]M \rightarrow 0.
\end{align*} 
is a graded minimal free resolution of $\gr^L[\bfn_0]M$ over $S$.

Moreover, by \cite[Proposition~4.2, 4.3, Corollary 4.1]{OT01}, the free resolution $\bfF$ also induces a graded minimal free resolution of a $D^\h$-module $M^\h$ 
\begin{align*}
\mathbf{F^\h}:\quad \cdots \xrightarrow[]{\psi_2} (D^\h)^{r_1}[\bfn_1]\xrightarrow[]{\psi_1} (D^\h)[\bfn_0] \xrightarrow[]{\psi_0} M^\h \rightarrow 0
\end{align*} such that the dehomogenization $M^\h |_{h=1} =M$ and $\bfF^\h |_{h=1} = \bfF$, where the maps $\psi_i$ are obtained by homogenizing the maps $\varphi_i$ for each $i \in \NN$.

Note that, $i >0$ and $j\in \ZZ$, we still have
\(\psi_i \left( L_j [\bfn_i](D^\h)^{r_i}\right) \subseteq  L_j [\bfn_{i-1}](D^\h)^{r_{i-1}},\) and that the free resolution $\bfF^\h$ also induced a graded free resolution over $S[h]$
\begin{align*}
\gr^L\mathbf{F^\h}:\quad \cdots \xrightarrow[]{\bar{\psi_2}} \gr^L[\bfn_1](D^\h)^{r_1}\xrightarrow[]{\bar{\varphi_1}} \gr^L[\bfn_0](D^\h)^{r_0} \xrightarrow[]{\bar{\varphi_0}} \gr^L[\bfn_0]M^\h \rightarrow 0.
\end{align*} 
On the other hand, the maps $\bar\varphi_i$ and $\bar{\psi_i}$ are obtained by the $L$-initial forms of entries in $\varphi_i$ and $\psi_i$. Since, for $P \in D$,  the $L$-initial form $\ini_{\oneone} P$ in $S$ coincides with the $L$-initial form of the homogenization of $P$ in $S[h]$,  we have the following commutative diagram
\[
\begin{CD}
\gr^L[\bfn_i](D^\h)^{r_i}     @>\overline{\psi_i}>> \gr^L[\bfn_{i-1}](D^\h)^{r_{i-1}} \\
@VVV        @VVV\\
k[h] \otimes \gr^L[\bfn_{i}]D^{r_{i}}   @> \id \otimes \overline{\varphi_i}>>  k[h] \otimes\gr^L[\bfn_{i-1}]D^{r_{i-1}}.
\end{CD}
\]

As a consequence, we see that $\gr^L[\bfn_0]M^\h \cong k[h] \otimes\gr^L[\bfn_0]M $ and that the free resolution 
 \[\gr^L \mathbf{F}^\h \cong k[h] \otimes \gr^L \mathbf{F} \] is a graded minimal free resolution of $\gr^L[\bfn_0]M^\h$ over $S[h]$.
 We therefore obtain the following result.
\begin{Proposition}\label{prop:GB_Betti}
    With the above notation, the graded Betti numbers of  $M^\h$, $\gr^{\oneone}M$, and $\gr^{\oneone} M^\h$ coincide for any finite $D$-module $M$.
\end{Proposition}

\section{Cohen--Macaulay $D$-modules and $D^\h$-modules}\label{sec:CMDmod}

In this section, we again consider the $L=\oneone$ filtration of the $n^\mathrm{th}$ Weyl algebra $D$ and the homogenized Weyl algebra $D^\h$ over a field $k$ of characteristic $0$. Recall that  the $k$-algebras $D^\h,S =\gr^L D$, and $S[h] = \gr^L D^\h$ are $\NN$-graded defined by the total degree of elements.

\subsection{}
Write $A:=D^\h$. Since the associated grade algebra $\gr^L A=S[h]$ is a polynomial ring of $2n+1$ variables over~$k$, ring $A$ is left and right Noetherian. So we can apply the results in section~\ref{sec:homalg}. In particular, for a nonzero finite left $A$-module $M$, we have
\[ 0 \le j(M) \le \pd(M) \le \gldim(A) \le 2n+1 \quad \text{and}\quad d(M)+j(M) = 2n+1.\]

On the other hand, since $A$ is an Auslander-regular $\NN$-graded $k$-algebra (Remark~\ref{rem:A-regular}), by \cite[Theorem~6.3]{Lev92} we have
\begin{equation*}
    \gExt_A^i(k,A) = \begin{cases}
    0 ,\quad i \ne 2n+1\\
    k ,\quad i=2n+1.
\end{cases}
\end{equation*}
As a consequence, we have $\depth(A) = \gldim(A) = 2n+1$ and the $\chi^\circ$ condition in section~\ref{sec:noncom AB} is satisfied. Therefore, by Theorem~\ref{thm:noncomAB} we have
\[\pd(M)+\depth(M) = \depth(A)=2n+1\] 
for any finite graded left $A$-module $M$.

Let us make the following definition.
\begin{Definition}
    Let $M$ be a nonzero finite graded left $D^\h$-module, we say that $M$ is {\it Cohen--Macaulay} if $\depth(M) = d(M)$ or equivalently if  $j(M)=\pd(M)$.
\end{Definition}

\subsection{}
Let $M$ be a nonzero finite left $D$-module. By Proposition~\ref{prop:GB_Betti}, there exists a graded $D^\h$-module $M^\h$ such that
the graded Betti numbers of $M^\h$, $\gr^L M$, and $\gr^L M^\h$ coincide. We abbreviate $\gr^L M$ by $\gr M$, and $\gr^L M^\h$  by $\gr M^\h$.
As a consequence, we have 
\begin{align*}
    &\pd(M^\h) = \pd(\gr M) = \pd(\gr M^\h) \le 2n, \text{ and} \\
     n+1 \le \,& d(M^\h) = \dim (\gr M^\h) = \dim(\gr M)+1 \le 2n+1.
\end{align*}
Note also that 
\begin{align*}
     d(M^\h)+j(M^\h) &= \pd M^\h) +\depth(M^\h) = 2n+1, \\
     \dim(\gr M^\h) + j(\gr M^\h) 
     &=\pd(\gr M^\h) +\depth(\gr M^\h) = 2n+1, \text{ and} \\
      \dim(\gr M) +j( \gr M) &=  \pd(\gr M) +\depth(\gr M) = 2n.
\end{align*}

Therefore, we deduce that $M^\h$ is Cohen--Macaulay if and only if $\gr M^\h$ is Cohen--Macaulay if and only if $\gr M$ is Cohen--Macaulay.  This justifies the following definition.

\begin{Definition}\label{def:CMD}
    A finite left $D$-module $M$ is {\it Cohen--Macaulay} if $\gr^{\oneone} M$ is Cohen--Macaulay.
\end{Definition}

When $M = D/I$ is a cyclic left $D$-module with $I$ a left $D$-ideal, we have $M^\h = A/ I^\h$, $\gr M^\h = S[h]/ \ini_{(\mathbf{1},\mathbf{1})} I^\h$, and $\gr M = S /\ini_{(\mathbf{1},\mathbf{1})} I$. The $D$-module $D/I$ is Cohen--Macaulay if by definition the $S$-module $S/ \ini_{\oneone} I$ is Cohen--Macaulay.

\begin{Remark}
Note that $\ini_{\oneone} I$ is a homogeneous ideal in $S$. 
Note that the ground field $k$ is infinite. Given a homogeneous ideal $J$ in $S$, it is well-known that if $S/\ini_{\prec}(J)$ is Cohen-Macaulay with respect to a term order $\prec$ then  $S/J$ is Cohen-Macaulay, see for example, \cite[Proposition 25.4]{PeevaBook}. The converse is not true in general. But if $\ini_{\prec}(J)$ is generated by square-free monomials, then  $S/J$ is Cohen-Macaulay implies $S/\ini_{\prec}(J)$ is Cohen-Macaulay by \cite[Corollary 2.7]{CV}.
\end{Remark}

\begin{Remark} 
    If a $D$-module $M$ is $L$-holonomic as in Remark~\ref{rem:1}(3), then $\dim(\gr M) =n$. In this case, the module $M$ is Cohen--Macaulay if and only if $j(\gr M) = \pd(\gr M) = \depth(\gr M) = n$. 
\end{Remark}

\section{Cohen--Macaulay GKZ systems}
In this section, we work over an algebraically closed field $\CC$ of characteristic $0$.

Let $A=(a_{ij})$ be a $d \times n$ integer matrix of rank $d$ whose columns $\bfa_1, \dots, \bfa_n \in \ZZ_d$ are nonzero. We assume that the semigroup $\NN A$ is pointed with $\ZZ A = \ZZ^d$. Consider the toric ideal $$I_A = \langle \partial^{u} - \partial^{v} \mid Au = Av\rangle$$   associated to $A$ in the polynomial ring $R= \CC[\partial_1, \dots, \partial_n]$ and the toric algebra $\CC[\NN A] \cong R/I_A$.

For $\beta \in \CC^d$, consider the GKZ-hypergeometric ideal
$H_A(\beta):= D\langle I_A, E - \beta \rangle$ in the Weyl algebra $D= \CC[x_1, \dots, x_n]\langle \partial_1, \dots, \partial_n \rangle$ where the Euler vector fields $E= (E_1, \dots, E_d)$ of $A$ are defined by $E_i:= \sum_j a_{ij}x_j \partial_j$ for $i = 1, \dots, d$. The $A$-hypergeometric system with parameter vector $\beta$ is the $D$-module $M_A(\beta) := D/ H_A(\beta)$.

Note that for any weight vector $(u,v) \in \ZZ^{2n}$ with $u_i+v_i>0$, the initial form $\ini_{(u,v)} (E_i-\beta_i)$ is equal to the sum $\sum_j a_{ij}x_j \partial_j$ in the polynomial ring $\CC[x,\d]$. We will sometimes interpret  $E_i$ as a polynomial in $\CC[x,\d]$ or identify $E_i$ with $\ini_{(u,v)}(E_i -\beta_i)$ for convenience.

Let $\ini I_A$ be the initial ideal of $I_A$ with respect to the total degree on $R=\CC[\d]$. The purpose of this section is to prove the following theorem. Recall that $S= \gr^{\oneone}D =\CC[x, \d]$.
\begin{Theorem}\label{thm:CM_GKZ}
    If $R/ \ini I_A$ is a Cohen--Macaulay $R$-module, then $\gr^{\oneone} M_A(\beta) = S/ \ini_{\oneone}H_A(\beta)$ is a Cohen--Macaulay $S$-module, that is, $M_A(\beta)$ is a Cohen--Macaulay $D$-module. 
\end{Theorem}

\begin{Remark} 
\begin{enumerate}
    \item As a special case of Theorem~\ref{thm:CM_GKZ}, if $I_A$ is a homogeneous ideal in $R$ so that $R/I_A$ is Cohen--Macaulay, then $M_A(\beta)$ is Cohen--Macaulay.
    \item It is known that if the Gr\"obner deformation $\CC[\d]/\ini I_A$ is Cohen--Macaulay, then so is $\CC[\d]/I_A$ \cite[Proposition~1.6.2]{BCRV}. However, the converse of this statement is not true in general (see \Cref{CMExample}).
    \item The Cohen--Macaulay condition on the toric ring $\CC[\d]/I_A$ is equivalent to the absence of rank-jumps in the GKZ hypergeometric system $H_A(\beta)$ \cite[Corollary~9.2]{MMW05}. It is also known in \cite[Corollary~4.13]{SW08} that if $\beta$ is not a rank-jumping parameter, the $L$-characteristic cycles of $M_A(\beta)$ is independent of $\beta$. In this case, the Cohen--Macaulay condition of $M_A(\beta)$ does not dependent on $\beta$.
    \item The converse statement of Theorem~\ref{thm:CM_GKZ} does not hold in general. See \Cref{CMExample}.
\end{enumerate}

\end{Remark}

\subsection{}
 Let $L = (L_{\d_1}, \dots, L_{\d_n}) \in \QQ^n$ be a weight vector on $R$. This induces an increasing filtration $L$ of $\CC$-vector space on $R$ by $[\d^u \in L_i R] \Leftrightarrow [L \cdot u \le i]$.
For $f =\sum c_u \d^u \in L_lR \setminus L_{l-1}R$, we define the $L$-degree  of $f$ as $\deg^L (f):=l$ and the $L$-initial form of   $f$  as $\ini_L(f) = \sum_{L\cdot u =l} c_u \d^u$. Also, for an ideal $I$ of $R$, we define the initial ideal $\ini_L I = \langle \ini_L(f) \mid f \in I \rangle$. The associated graded ring $\gr^L R$ is isomorphic to $R$. By abuse of notation, we identify the associated graded module $\gr^L (R/I)$ with $R / \ini_L I$.

\subsection{}
The geometry of $\Spec(R/ \ini_LI)$ is characterized by the $(A,L)$-umbrella $\Phi^L_A$ introduced in \cite[Definition~2.7]{SW08}, which is the set of faces of $\conv_H (\{0, \bfa_1^L, \dots,\bfa_n^L\})$ that do not contain zero where $\bfa_i^L = \bfa_i /L_{\d_i}$ for $i=1,\dots, n$. Here we embed $\{0, \bfa_1^L, \dots,\bfa_n^L\} \subset \QQ^d$ into $\PP^d_\QQ$ and consider the convex hull of 
$\{0, \bfa_1^L, \dots,\bfa_n^L\}$ in $\PP^d_\QQ$ relative to a hyperplane $H$ that separates $\{0\}$ and $\{\bfa_1^L, \dots,\bfa_n^L\}$.

For $\tau \in \Phi_A^L$, we identify it with $\{j \mid \bfa_j^L \in \tau\}$, or with $\{\bfa_j \mid \bfa_j^L \in \tau\}$, or with the corresponding submatrix of $A$. By $\Phi^{L,d-1}_A$, we denote the subset of faces of dimension $d-1$. Denote by $R_\tau = \CC[\d_\tau]$ the polynomial subring of $R$ generated by the variables $\d_i$, $i \in \tau$. Consider the toric ideal $I_\tau$ in $R_\tau$ and the ideal $J_\tau:= \langle \d_i \mid i \notin \tau \rangle$ in $R$. On the other hand, the $d$-torus $\TT:= (\CC^*)^d$ with coordinates $t:= (t_1, \dots, t_d)$ acts on $\Spec R$ by $t\cdot \xi:= (t^{\bfa_1} \xi_1, \dots, t^{\bfa_n} \xi_n)$, which induced a $\ZZ^d$-grading on $R$ by $\deg(\d_i)= \bfa_i$. Define $\mathbf{1}^\tau_A \in \Spec R$ by $(\mathbf{1}^\tau_A)_i := 1$ if $\bfa_i \in \tau$ and $(\mathbf{1}^\tau_A)_i := 0$ if $\bfa_i \notin \tau$. Denote by $O^\tau_A$ the $\TT$-orbit of $\mathbf{1}^\tau_A$ and by $\overline{O^\tau_A}$ its Zariski closure.

\begin{Theorem}{\cite[Theorem~2.14]{SW08}}\label{thm:Var(inI_A)}
    The set of $\ZZ^d$-graded prime ideal of $R$ containing $\ini_L I_A$ equals $\{ RI_\tau +J_\tau \mid \tau  \in \Phi^L_A\}$, and we have
    \[\Spec(R/ \ini_L I_A) = \Var(\ini_L I_A) = \bigcup_{\tau \in \Phi^{L,d-1}_A} \overline{O^\tau_A} = \bigsqcup_{\tau \in \Phi^L_A} O^\tau_A. \]
\end{Theorem}

Notice that Theorem~\ref{thm:Var(inI_A)} holds when replacing $\CC$ by any algebraically closed field $K$ of characteristic $0$. Let $K=\overline{\CC(x)}$ be the algebraic closure of $\CC(x)=\CC(x_1, \dots, x_n)$. We have the following lemma.
\begin{Lemma}\label{lem:lsop}
The Euler operators $E_1, \dots, E_d$ form a linear system of parameters of $\CC(x)[\d]/ \ini_L(I_A)$.
\end{Lemma}
\begin{proof}
    It follows from Theorem~\ref{thm:Var(inI_A)} that each associated primes of $\ini_L I_A$  is of codimension $d$, so we have 
 $ \dim \left( R/ \ini_L I_A\right)=d.$
    Moreover, since $\CC(x)[\d]/ \ini_L(I_A) = \left(R/ \ini_L I_A \right) \otimes \CC(x)$, we also have $\dim \left(\CC(x)[\d]/ \ini_L(I_A) \right) = d$. So it suffices to show that $\CC(x)[\d]/ \left(\ini_L(I_A)+ \langle E \rangle \right)$ is supported at the origin, namely  \(\Var_{\CC(x)}\left( \ini_L I_A +\langle E \rangle \right) = \{0\}.\)

    Apply Theorem~\ref{thm:Var(inI_A)} to $K[\d] / \ini_L I_A$, we have $\Var_K(\ini_L I_A) = \bigsqcup_{\tau \in \Phi^L_A} O^\tau_A$.
    Let $\tau \in \Phi^L_A$ be such that $E_1,\dots,E_d$ all vanish at $\mathbf{1}^\tau_A \in \Var_K(\ini_L I_A)$. If $\tau \ne \emptyset$, rearrange $\bfa_1, \dots, \bfa_n$ so that $\bfa_1 \in \tau$ and that $\bfa_1, \dots, \bfa_d$ are $\QQ$-linearly independent. Let $\sigma$ be the submatrix of $A$ whose columns are $\bfa_1, \dots, \bfa_n$. Denote by $(s_1,\dots,s_d)$ the first row of the invertible matrix $\sigma^{-1}$. Then the operator \[\sum_{i=1}^d s_i E_i \in \left( x_1\d_1 +\sum_{j=d+1}^n \QQ x_j \d_j \right)\] also vanishes at $\mathbf{1}^\tau_A$. From this we deduce that $x_1 \in \sum_{i \in \Lambda} \QQ x_i$ for some nonempty $\Lambda \subseteq \{d+1,\dots, n\}$, which is not possible.
\end{proof}

The following lemma is motivated by \cite[Lemma~4.3.7]{SST_book}.
\begin{Lemma}\label{lem:Euler_rs}
    For a matrix $A$ and a weight vector $L$ on $R=\CC[\d]$, the following conditions are equivalent :
\begin{enumerate}
    \item $\CC[\d]/\ini_L(I_A)$ is Cohen--Macaulay;
    \item $\CC(x)[\d]/\ini_L(I_A)$ is Cohen--Macaulay;
    \item $\CC[x][\d]/\ini_L(I_A)$ is Cohen--Macaulay;
    \item The Euler operators $E_1, \dots, E_d$ form a regular sequence on $\CC[x][\d]/\ini_L(I_A)$.
    \item The Euler operators $E_1, \dots, E_d$ form a regular sequence on $\CC(x)[\d]/\ini_L(I_A)$.
\end{enumerate}
\end{Lemma}
\begin{proof}
    Since $\CC(x)[\d]/ \ini_L(I_A) = \left(R/ \ini_L I_A \right) \otimes \CC(x)$ and $\CC[x][\d]/ \ini_L(I_A) = \left(R/ \ini_L I_A \right) \otimes \CC[x]$,
    the equivalence of (1), (2), and (3) follows from \cite[Theorem~2.1.9, 2.1.10]{BH93}.

    Let $P_\d$ be the prime ideal in $\CC[x][\d]/ \ini_L I_A$ generated by the images of $\d_1, \dots, \d_n$. Notice that the multiplicative closed set of elements not in $P_\d$ is equal to $\{ f + \ini_L I_A \mid f \in \CC[x] \setminus \{0\}\}$. So  $\CC(x)[\d]/\ini_L(I_A)$ is isomorphic to the localization $\CC[x][\d]/\ini_L(I_A)$ at $P_\d$ , namely
    \begin{equation}\label{eq:loc_ini}
        \CC(x)[\d]/\ini_L(I_A)=\left(\CC[x][\d]/\ini_L(I_A)\right)_{P_\d}. 
    \end{equation}
    Since the regular sequence is preserved by localization, we have (4) implies (5).

    Since $\dim \left(\CC(x)[\d]/\ini_L(I_A)\right)=d$, by \cite[Theorem~1.2.5]{BH93} the existence of a length $d$ regular sequence implies that $\depth  \left(\CC(x)[\d]/\ini_L(I_A)\right) =d$. This shows that (5) implies (2).

    We conclude the proof by showing that (3) implies (4). Let $J_E \subset P_\d$ be the ideal in $\CC[x][\d]/ \ini_L I_A$ generated by the image of $E_1, \dots, E_d$.
    By Lemma~\ref{lem:lsop} and \eqref{eq:loc_ini}, we have
    \[\height J_E = \height P_\d  = \height (J_E)_{P_\d} = \dim \left(\CC(x)[\d]/\ini_L(I_A)\right)=d. \] Since $\CC[x][\d]/\ini_L(I_A)$ is Cohen--Macaulay by (3), it follows from \cite[Theorem~2.1.6]{BH93} that the ideal $J_E$ is unmixed and that $E_1, \dots, E_d$ form a regular sequence on $\CC[x][\d]/\ini_L(I_A)$.
\end{proof}
\vspace*{-0.5cm}
We are now ready to prove Theorem~\ref{thm:CM_GKZ}.
\begin{proof}[Proof of Theorem~\ref{thm:CM_GKZ}]
    Recall that $R=\CC[\d]$ and $S= \gr^{\oneone}D =\CC[x, \d]$. For all one vector $\mathbf{1} \in \ZZ^n$, notice that $\ini I_A = \ini_{(\mathbf{1})} I_A$ in $R$ and that $\ini_{\oneone}I_A = S(\ini I_A) $ in $S$. 

    Since $R/ \ini I_A$ is Cohen--Macaulay, by Lemma~\ref{lem:Euler_rs} the initial forms $\ini_{\oneone}(E-\beta)= E$ is a regular sequence on $S/\ini_{\oneone} I_A$. So the assumption of \cite[Theorem~4.3.5]{SST_book} is satisfied and we have 
\[\ini_{(1,1)} H_A(\beta) = \ini_{(1,1)} I_A + \langle E \rangle.\]
Since $S/\ini_{\oneone} I_A = \left( R/ \ini I_A \right) \otimes \CC[x]$ is Cohen--Macaulay, it follows from \cite[Theorem~2.1.3]{BH93} that $S/ \ini_{\oneone}H_A(\beta)$ is Cohen--Macaulay.
\end{proof}

We conclude this article with examples illustrating the relations between the Cohen--Macaulay condition on the three modules $\CC[\NN A] = R/I_A$, $R/\ini I_A=R/\ini_{(\mathbf{1})} I_A$, and 
\[\gr^{\oneone}M_A(\mathbf{0})=S/\ini_{\oneone} H_A(\mathbf{0})= S/\ini_{\oneone} D \langle I_A, E \rangle.\]
\begin{Example} \label{CMExample}
The following table includes 5 different examples when $\dim \CC[\NN A]=2$.
\vskip 20pt
\begin{center}
    \begin{tabular}{|c|c|c|c|}
\hline 
 & $\CC[\NN A]$ & $R/\ini I_A$ & $S/\ini_{\oneone} D \langle I_A, E_1, E_2 \rangle $\tabularnewline
\hline 
\hline 
$A=\left(\begin{array}{cccc}
0 & 1 & 2 & 2\\
2 & 1 & 1 & 0
\end{array}\right)$ & CM & CM & CM\tabularnewline
\hline 
$A=\left(\begin{array}{cccc}
0 & 1 & 2 & 3\\
1 & 1 & 0 & 0
\end{array}\right)$ & not CM & not CM & CM\tabularnewline
\hline 
$A=\left(\begin{array}{cccc}
0 & 1 & 2 & 3\\
1 & 2 & 0 & 0
\end{array}\right)$ & not CM & not CM & not CM\tabularnewline
\hline 
$A=\left(\begin{array}{cccc}
0 & 1 & 2 & 3\\
1 & 2 & 2 & 0
\end{array}\right)$ & CM & not CM & not CM\tabularnewline
\hline 
$A=\left(\begin{array}{ccccc}
1 & 2 & 3 & 4 & 5\\
3 & 1 & 3 & 2 & 0
\end{array}\right)$ & CM & not CM & CM\tabularnewline
\hline 
\end{tabular}
\end{center}
\vskip 20pt

Here, since $\dim \CC[\NN A]=2$, it follows from \cite[Corollary~6.2.6]{BH93}  that $\CC[\NN A]$ is Cohen--Macaulay if and only if $H^2_\frakm(\CC[\NN A])$ is the only nonvanishing local cohomology of $\CC[\NN A]$. Moreover, the local cohomology of the semigroup ring $\CC[\NN A]$ can be computed by the Ishida complex \cite[Theorem~6.2.5]{BH93}. More concretely, let $\bfa_1$ and $\bfa_2$ be the columns of $A$ that lie on the two rays of the cone $\QQ_{\ge 0}A$. The Cohen--Macaulay condition of $\CC[\NN A]$ is equivalent to the condition $(\NN A -\NN \bfa_1) \cap (\NN A -\NN \bfa_1) = \NN A$. 

On the other hand, we use the Auslander-Buchsbaum
formula and the projective dimension of $R/\ini I_A$
or $S/\ini_{\oneone}  \langle I_A, E_1, E_2 \rangle $
to test their Cohen-Macaulay condition with the help of Macaulay2 \cite{M2} and the Dmodules package written by Anton Leykin and Harrison Tsai. 
\end{Example}

\bibliography{Ref}
\end{document}